# MUTUAL INTERPRETABILITY OF ROBINSON ARITHMETIC AND ADJUNCTIVE SET THEORY WITH EXTENSIONALITY

Zlatan Damnjanovic

University of Southern California

[Abstract: An elementary theory of concatenation, $QT^+$, is introduced and used to establish mutual interpretability of Robinson arithmetic, Minimal Predicative Set Theory, quantifier-free part of Kirby's finitary set theory, and Adjunctive Set Theory, with or without extensionality.]

The theory Q, also known as Robinson arithmetic (described in sec. 12 below), is often singled out as of special interest in foundational arguments, for different reasons. In textbook presentations of Gödel's incompleteness theorems and related undecidability results, Q is introduced as a minimal formal deductive framework allowing for arithmetical representation of formalized syntax and definition of basic concepts of recursion theory, furnishing an example of a finitely axiomatizable yet essentially undecidable theory. (See, e.g., [16].) Mathematically, it appears to be extremely weak: even though addition and multiplication are present, their most basic properties, such as associativity and commutativity, cannot be formally proved in it. This is hardly surprising, considering that induction is not included among the axioms. Nonetheless, the notion of a theory relatively interpretable in Q – in the sense explained, e.g., in [16] – proved to be surprisingly rich and of great independent interest in connection with a philosophically motivated neo-formalist program put forward by Edward Nelson. (See [12]. Nelson's position is sometimes also described as strict predicativism.) By a theorem of Wilkie, the absence of induction among the axioms of Q does not preclude the possibility of interpreting in Q a fragment of Peano Arithmetic, I$\Sigma_0$, where induction is restricted to formulas with bounded quantifiers; this allows recovery (via the interpretation) of all the usual arithmetical laws. Additionally, an impressive amount of non-trivial mathematics can be reconstructed in theories interpretable in Q, including (first-order) Euclidean geometry, elementary theory of the real closed fields (i.e., first-order theory of real numbers) as well as basic "feasible analysis" formalizing elementary properties of real numbers and continuous functions. (See [4] for details.) It is frequently pointed out that Q is a minimal element in the well-ordered hierarchy of interpretability of "natural" mathematical theories.

It was Tarski who first noted that, as regards self-referential constructions at the heart of meta-mathematical arguments for incompleteness, the procedure of arithmetization by means of which the syntax of formal theories is coded up by numbers amounts to an unnecessary detour. In his seminal work on the concept of truth of formalized languages Tarski introduced a theory of concatenated strings to demonstrate this point. This idea was further developed by Quine [13]. More recently, Grzegorczyk has suggested that a theory of concatenated "texts" would form a natural framework for the study of incompleteness phenomena and, more generally, computation, and for this purpose he introduced a weak theory of concatenation, TC, and proved its undecidability [6]. Work of several authors has since revealed an intimate connection between TC and Q: the two theories turn out to be mutually interpretable, as established independently by Visser and Sterken [19], Ganea [5], and Švejdar [15]. For a thorough discussion of interpretability in connection with these theories, see [19].

We approach these issues from a somewhat different angle. Is there a comparably "minimal" theory of sets that is deductively equipollent, in the sense of mutual interpretability, to Q? This problem has a long and interesting history. A hint of a positive answer was provided by Szmielew and Tarski in 1950, who announced in [14] interpretability of Q in a very weak fragment of set theory called Adjunctive Set Theory with Extensionality, AST+EXT (this theory is described in sec. 5 below). The result was restated in [16], p.34, but no proof was published. Collins and Halpern did produce a proof in 1970, said to be approved of by Tarski (see [2]). In 1990 an interpretation of Q in AST dispensing with Extensionality was outlined in Appendix III of [10], by Mycielski, Pudlák and Stern. In 1994 Montagna and Mancini, who seem not to have been aware of [2] and [10], showed that Q is interpretable in an extension of AST they call Minimal Predicative Set Theory, N. The theory N was said to be suggested to them by Nelson, with a constant for the empty set, a binary relation symbol ∈ and a binary function symbol for set adjunction x∪{y}, but without Extensionality (see [11]; this theory is described in sec. 14 below). A new proof of interpretability of Q in AST was given in Visser [18], and also, in a somewhat different setting, in Burgess [1].

The reverse problem of interpreting AST in Q appears more vexing because of the paucity of resources available in Q for set construction. For discussion of the important notion of sequential theory broadly relevant in this connection see [8], [10], and especially the works of Visser [17] and [19]. One path would be to try to think of sets in terms of numbers and then model set-building operations – in this case adjunction – arithmetically, modulo a sufficiently elementary coding, relying on bounded induction in a suitable extension of $I\Sigma_0$ known to be interpretable in Q, such as, e.g., $I\Sigma_0+\Omega_1$, to set up the needed interpretation. Indeed, along these lines, one can obtain an interpretation of a version of Adjunctive Set Theory (without Extensionality) in Q. (This was essentially accomplished by Nelson in [12]). We follow a different tack. Building on some ideas of Quine, we propose to think of sets in terms of strings, and of adjunction as concatenation or juxtaposition of certain kinds of strings. We introduce a weak theory of concatenation, QT, different from Grzegorczyk's TC, to intermediate between Adjunctive Set Theory and Q.

Our approach will turn out to have the unexpected philosophical benefit of allowing us to understand both arithmetic and set theory as just forms of concatenation theory. The theory QT, or rather its definitional extension QT$^+$, will serve as an interpretive framework for set adjunction, while QT$^+$ can in turn be interpreted in I$\Sigma_0$. (It can be shown independently that Q is directly interpretable in QT$^+$.) Here we confront two principal challenges: in the absence of induction in QT, to establish the relevant facts about the coding of sets of strings by strings needed to carry out the interpretation, and, secondly, to ensure that the resulting interpretation validates the extensionality axiom for sets. The simultaneous solution of these two problems yields our interpretation of AST+EXT in QT$^+$. Since I$\Sigma_0$ is interpretable in Q, and, as shown in sec. 11 below, QT$^+$ is interpretable in I$\Sigma_0$, this suffices to establish mutual interpretability of AST+EXT and Q.

Some authors prefer to formalize set adjunction using an adjunction operator represented by a primitive binary function symbol for x∪{y}, rather than relying on the logical apparatus of the first-order language of set theory $\mathcal{L} = \{\in\}$. One such formulation, called PS (for "Peano Set Theory"), proposed by Kirby in [9], is intended to serve as an axiomatization of the first-order theory of hereditarily finite sets. Analogously to how first-order Peano Arithmetic can be introduced as an extension of Q by means of the induction schema over the language {0, ′,+,·} of arithmetic, PS is presented as a system, $PS_0$, consisting of a handful of quantifier-free axioms extended by an induction schema expressed in terms of the adjunction operator. The axioms of the system N of Minimal Predicative Set Theory of Montagna and Mancini are included among those of $PS_0$, which is described below in sec. 14.

We establish that Robinson arithmetic, Minimal Predicative Set Theory, Kirby's quantifier-free finitary set theory, and Adjunctive Set Theory, with or without extensionality, are all relatively interpretable in each other, each being mutually interpretable with concatenation theory $QT^+$. Thus, fundamentally, the most basic arithmetic and simplest set theory turn out to be variants of one and the same theory.

Our arguments requires construction of formulas with certain special properties, and tedious formal verifications of those properties. While we provide the reader with sufficient indications how these constructions can be carried out, the complete details are given in [3].

The author is grateful to the anonymous referee for very helpful comments.

## 1. A Theory of Concatenation

We consider a first-order theory with identity and a single binary function symbol *. Informally, we let the variables range over nonempty strings of a's and b's – or 0's and 1's – and let x*y be the string that consists of the digits of the string x followed by the digits of string y, subject to the following conditions:

(QT1)  (x*y)*z = x*(y*z)

(QT2)  ¬(x*y=a) & ¬(x*y=b)

(QT3)  (x*a=y*a → x=y) & (x*b=y*b → x=y) &

            & (a*x=a*y → x=y) & (b*x=b*y → x=y)

(QT4)  ¬(a*x=b*y) & ¬(x*a=y*b)

(QT5)  x=a v x=b v (∃y(a*y=x v b*y=x) & ∃z(z*a=x v z*b=x))

It is convenient to have a function symbol for a successor operation on strings:

(QT6)    $Sx=y \leftrightarrow ((x=a\ \&\ y=b) \vee (\neg x=a\ \&\ x*b=y))$.

Of course, we may also think of a single letter a appended to x, x*a, as a successor of the string x.

Since the last axiom is basically a definition, adding it to the other five results in an inessential (i.e. conservative) extension.

We call the resulting theory $QT^+$.

First, let's introduce some obvious abbreviations:

   $xBy \equiv \exists z\ x*z=y$   and   $xEy \equiv \exists z\ z*x=y$.

Also, let    $x \subseteq_p y \equiv x=y \vee xBy \vee xEy \vee \exists y_1 \exists y_2\ y=y_1*(x*y_2)$.

Often, we write xy for x*y and omit parentheses in x(yz) and (xy)z on account of (QT1).

For philosophical reasons, and at a slight cost of increase in complexity of definitions such as those just given, we do not include the empty string. It can be shown, however, that concatenation theory with the empty string and one without are mutually interpretable (see [7]).

Let    $xRy \equiv (x=a\ \&\ \neg y=a) \vee xBy$.

Then we can prove (see [3], pp.13-17):

$$QT^+ \vdash xR(Sx),$$

$$QT^+ \vdash Sx=Sy \rightarrow x=y,$$

$$QT^+ \vdash \neg xRa \,\&\, (\neg x=a \rightarrow aRx)$$

$$QT^+ \vdash xRy \,\&\, yRz \rightarrow xRz$$

$$QT^+ \vdash xR(Sy) \leftrightarrow xRy \lor x=y$$

These elementary facts tell us that, provably in $QT^+$,

$xRy \lor x=y$ is a discrete preordering of strings.

Still, our theory looks exceedingly weak.

## 2. Tractable strings

Consider the following problem. We know that neither $a$ nor $b$ are their own initial segments:

$$QT^+ \vdash \neg aRa \quad \text{and} \quad QT^+ \vdash \neg bRb.$$

But we <u>don't know whether our theory proves that no string</u> is an initial segment of itself: $\qquad QT^+ \vdash^? \forall x \,\neg xRx\,.$

Let's look at this more carefully.

Consider the property $\quad I_0(y) \equiv \forall x\,(xRy \lor x=y \rightarrow \neg xRx)$.

Both a and b have this property:

$$QT^+ \vdash I_0(a) \text{ and } QT^+ \vdash I_0(b).$$

By definition, <u>no string in $I_0$</u> is its own initial segment:

(!) $\qquad\qquad QT^+ \vdash I_0(x) \rightarrow \neg xRx.$

Let us call $I_0$ strings <u>tractable</u>. So a and b are tractable.

It may be that not all strings are tractable.

But we are going to be working with those that are.

To begin, together with transitivity of the relation R, we have that (!) delivers anti-symmetry for tractable strings ([3], (2.2)):

$$QT^+ \vdash I_0(x) \rightarrow \neg(xRy \,\&\, yRx).$$

So our theory proves that tractable strings form a discrete partial ordering under R with a as the least element.

Write $x<y$ for $I_0(x)\,\&\,I_0(y)\,\&\,xRy$.

As usual, $x \leq y$ stands for $x<y \lor x=y$.

Let us summarize:

<u>Theorem 0.</u> (1)  $QT^+ \vdash \forall x \in I_0\ a \leq x$,

(2)  $QT^+ \vdash \forall x,y,z \in I_0\ (x \leq y\ \&\ y \leq z\ \rightarrow\ x \leq z)$,

(3)  $QT^+ \vdash \forall x \in I_0\ (x \leq Sx\ \&\ \neg x = Sx)$,

(4)  $QT^+ \vdash \forall x,y \in I_0\ (x \leq Sy\ \leftrightarrow\ x \leq y \lor x = Sy)$,

(5)  $QT^+ \vdash \forall x,y \in I_0\ (Sx = Sy \rightarrow x = y)$,

(6)  $QT^+ \vdash \forall x,y \in I_0\ (x \leq y\ \&\ y \leq x\ \rightarrow\ x = y)$.

We write "$\forall x \in K$" for "$\forall x\ (K(x) \rightarrow ...)$".  This relativization notation will prove to be convenient.

## 3. String concepts

Other than  a  and  b, which strings are tractable?

We don't know yet, for instance, whether, according to our theory, the tractable strings are closed under the successor operations on strings.

It turns out that they are ([3], pp.18-21):

   $QT^+ \vdash I_0(x) \rightarrow I_0(Sx)$   and   $QT^+ \vdash I_0(x) \rightarrow I_0(x*a)$.

Formulae with this property will be of special interest to us.

In general, we call a formula I(x) in the language of concatenation theory a string concept if, provably in $QT^+$, it holds of the letter a and hereditarily w.r. to both successor operations:

$QT^+ \vdash I(a)$,    $QT^+ \vdash \forall x\, (I(x) \rightarrow I(Sx))$  and  $QT^+ \vdash \forall x\, (I(x) \rightarrow I(x*a))$.

Note that if $I_1(x)$ and $I_2(x)$ are string concepts, so is their conjunction.

Now, $I_0(x)$ is a string concept.

Of course, so is  x=x !

But not knowing whether our theory proves that every string is tractable,

whether  $QT^+ \vdash \forall x\, \neg xRx$,  we didn't have the analogue of (!),  that no string to which the concept applies is its own initial segment, which gave us

(3) and (6) in Theorem 0.

In fact, $QT^+ \nvdash \forall x\, \neg xRx$. A countermodel exists. Consider, for example, a model $M^\infty$ of $QT^+$ with an infinite word $W \in M$ where

$$W = bb\ldots\ldots bbaabb\ldots\ldots bbaabb\ldots\ldots bb,$$

which begins and ends with an infinite sequence of b's and has midsegments …bbaabb… that begin and end with an infinite sequence of b's and form a countable dense linear ordering without endpoints. Then $X \in M^\infty$ where

$$X = abb\ldots .bbaabb\ldots \ldots bbaabb\ldots \ldots bba,$$

so that aWa=X. By a theorem of Cantor, there is a 1-1 order-preserving map from X onto any of its proper initial segments that end with ba. So $M^\infty \vDash X*X=X$ and thus $M^\infty \vDash XBX$. Therefore $QT^+ \nvdash \forall x \neg xBx$:

## 4. Simulating induction

Let's see how string concepts can help us cope with the apparent deductive weakness of our concatenation theory.

Our theory does not include an axiom or schema of induction.

Is there some way to rely on reasoning about string concepts to derive non-trivial universal conclusions about strings?

Let's consider an example.

Axiom (QT3) gives us left and right cancellation of <u>atoms</u>.

What about right cancellation for <u>all</u> strings

$$\forall z\ \forall x,y\ (x*z=y*z \rightarrow x=y)\ ?$$

Suppose we have a string concept I such that right cancellation holds for some string u ∈ I:

(hyp u) $\quad\quad\quad\quad \forall x,y \in I \ (x*u=y*u \rightarrow x=y)$

What about <u>its b-successor</u>, the string u*b? Can we cancel it?

Note that, we do have, for x, y ∈ I,

$x*(u*b)=y*(u*b) \Rightarrow (x*u)*b=(y*u)*b \Rightarrow x*u=y*u \Rightarrow x=y$

by (QT1), (QT3) and (hyp u). Likewise with a.

Thus

$QT^+ \vdash \forall x,y \in I \ (x*u=y*u \rightarrow x=y) \rightarrow \forall x,y \in I \ (x*(u*b)=y*(u*b) \rightarrow x=y)$,

and likewise for a-successors.

But this doesn't entitle us to conclude that

$\forall z \in I \ \forall x,y \in I \ (x*z=y*z \rightarrow x=y)$.

What if z happens to be an "infinite" string, one that cannot be obtained from an atom by repeatedly concatenating b or a? Nothing guarantees that such a string, if it exists, will in fact be cancellable. So we cannot claim, based on our theory, that all I-strings are right cancellable.

Suppose we contemplate a new string concept, J(z), using the very formula that we want to come out universally true:

($)     $J(z) \equiv I(z) \;\&\; \forall x,y \in I \;(x*z=y*z \to x=y).$

It is easy to see that the new predicate holds for a.

Likewise for b: $QT^+ \vdash J(b)$.

Now, assume $J(z)$ and let $x*(z*a)=y*(z*a)$. Then $(x*z)*a=(y*z)*a$ by (QT1). But then $x*z=y*z$ by (QT3), whence $x=y$ by hypothesis $J(z)$.

On the other hand, from $J(z)$ we have $I(z)$, and given that I is a string concept, $I(z*a)$.

Therefore,

$$QT^+ \vdash \forall z \;(J(z) \to J(z*a)).$$

Completely analogously, $QT^+ \vdash \forall z \;(J(z) \to J(z*b))$.

Here is the point. Obviously, that all J-strings are right cancellable,

($$)     $QT^+ \vdash \forall z \in J \;\forall x,y \in J \;(x*z=y*z \to x=y),$

follows from the definition ($).

But now we know that J is a string concept!



If we take I to be $I_0$, then we now know that, among the tractable strings, those with the right cancellation property include not only a and b, but are also, provably in our theory, closed under both successor operations. We established this by refining the given string concept I, then formulating the universal proposition for the resulting predicate J, and then proving that J is a string concept.

To summarize, faced with the prospect that the property of right cancellability may not hold universally for arbitrary I-strings z, we still managed to legitimately universalize by selecting those among I-strings that do have the property. We did this simply by writing out the appropriate condition that strengthens I to J, expressing the claim that the right cancellation can hold universally for I–strings in principle, albeit only for those – namely J-strings – that also satisfy the extra condition. Our theory guarantees that such strings are plentiful, because, as we saw in the above formal argument, $QT^+$ proves that J is a string concept if I is. If I-strings have been selected from $I_0$, i.e. tractable strings, this opens up the possibility of discovering other general properties of right cancellable tractable strings by strengthening these conditions by further requirements.



Some people may feel that this is cheating, that universalizing by strengthening is not really induction because the new string concept J is not the same as the old one, I, and not all I-strings may in fact be J-strings. If you feel that way, feel free to call it <u>quasi-induction</u> or simulated induction.

This procedure will allow us to simulate induction in our concatenation theory without explicitly assuming it.

## 5. Formula selection further illustrated

Our ability to exploit the extremely meager deductive resources of $QT^+$ has been potentially amplified by the knowledge that new universal statements, such as the one given in ($\$\$$), can be proved in $QT^+$ in a form relativized to a suitably selected string concept. We should note that the condition of right cancellability was, in logical terms, simple enough that we could count on $QT^+$ to deliver the needed conclusion that J is a string concept. To deal with other, more complex conditions it will turn out to be useful to know that the restricting condition J has some additional features.

Do we know that the strings in $I_0$ are closed under *? We don't. We need a string concept with that property:

CLOSURE UNDER CONCATENATION. There is a string concept $I_1$ such that
$$QT^+ \vdash \forall x\, \forall y\, (I_1(x)\, \&\, I_1(y) \rightarrow I_1(x*y))$$
where $QT \vdash \forall x\, (I_1(x) \rightarrow I_0(x))$.

Let $\quad I_1(x) \equiv I_0(x)\ \&\ \forall y(\, I_0(y) \rightarrow I_0(y*x))$.

We need to verify that $I_1(x)$ is indeed a string concept. First, that $QT^+ \vdash I_1(a)$: we have $QT^+ \vdash I_0(a)$. Suppose $I_0(y)$. Then $I_0(y*a)$ because $I_0$ is closed under $S_a$, provably in $QT^+$. So indeed $QT^+ \vdash I_1(a)$.

As for $QT^+ \vdash I_1(b)$, that follows from $QT^+ \vdash I_0(a)$ and the closure of $I_0$ under S. Next we show that $QT^+ \vdash \forall z\, (I_1(z) \rightarrow I_1(z*b))$.

So suppose $I_1(z)$. We want $I_1(z*b)$. We have $I_0(z)$ from the hypothesis $I_1(z)$, and so $I_0(z*b)$ since $I_0$ is a string concept. Assume that $I_0(y)$. From the hypothesis $I_1(z)$ it then follows that $I_0(y*z)$, and further that $I_0((y*z)*b)$. By (QT1), this means that $I_0(y*(z*b))$. So we have established that
$$\forall y(\, I_0(y) \rightarrow I_0(y*(z*b)),$$

which, along with the previously obtained $I_0(z*b)$, gives us $I_1(z*b)$ under the hypothesis $I_1(z)$, as required.

Similarly, $QT^+ \vdash \forall z\, (I_1(z) \to I_1(z*a))$.

This completes the argument that $I_1(x)$ is a string concept. But we also need to that $I_1$ is actually closed under the concatenation operation *, that is,

$$QT^+ \vdash \forall x\, \forall y\, (I_1(x)\ \&\ I_1(y) \to I_1(x*y)).$$

Assume $I_1(x)$ and $I_1(y)$, namely

(a) $\quad I_0(x)\ \&\ \forall z(I_0(z) \to I_0(z*x))$,

and (b) $\quad I_0(y)\ \&\ \forall z(I_0(z) \to I_0(z*y))$.

From $I_0(x)$ and (b) we obtain $I_0(x*y)$. Assume now $I_0(z)$. Then $I_0(z*x)$ by (a), and further $I_0((z*x)*y)$ by (b). But then $I_0(z*(x*y))$ by (QT1). So we have that $\quad I_0(x*y)\ \&\ \forall z(I_0(z) \to I_0(z*(x*y)))$, that is, $I_1(x*y)$. This is precisely what we needed to show.

Note that we have not used any property specific to $I_0$ as a string concept in the above argument. Say that a string concept I is <u>stronger than</u> $I_0$ if $QT^+ \vdash \forall x\, (I(x) \to I_0(x))$ and write $I \subseteq I_0$. We have in fact proved something more general: for any string concept I⊆I₀ there is a string concept J⊆I such that $\qquad$ QT⁺ ⊢ ∀x ∀y (J(x) & J(y) → J(x*y)).

Similarly, we can ensure that a given string concept is also downward closed with respect to the relation ≤ among tractable strings:

DOWNWARD CLOSURE UNDER ≤. Suppose J⊆I is a string concept where I⊆I₀. Then there is a string concept J≤⊆J such that

$\qquad$ QT⁺ ⊢ ∀x (J≤(x) & y≤x → J≤(y)).

Let J≤(x) ≡ ∀y≤x J(y).

We write ∀y≤x … for ∀y(y≤x → …).

That the formula J≤(x) has the required property is immediate from the definition and transitivity of ≤.

We have QT⁺ ⊢ J(a) by hypothesis, and QT ⊢ y≤a ↔ y=a. So QT⁺ ⊢ J≤(a).

Suppose M ⊨ J≤(x). Then M ⊨ ∀y≤x J(y). Suppose M ⊨ y≤Sx. Then

$\qquad$ y≤Sx ↔ y≤x ∨ y=Sx.

If y≤x, we have J(y) from the hypothesis J≤(x).

If y=Sx, then J(x) from the hypothesis J≤(x), whence J(Sx) from the principal hypothesis. Therefore, ∀y≤Sx J(y), that is, J≤(Sx). That J≤(x) is closed under $S_a$ is proved in the same fashion. This completes the argument that $J_\leq(x)$ is a string concept.

More generally, we have:

<u>DOWNWARD CLOSURE UNDER SUBSTRINGS.</u> For any string concept $I \subseteq I_0$ there is a string concept $J \subseteq I$ such that

$$QT^+ \vdash \forall x \in J \; \forall y \; (y \subseteq_p x \rightarrow J(y)).$$

Let $I^{\subseteq p}(x) \equiv I(x) \; \& \; \forall z \leq x \; \forall y \; (y \subseteq_p z \rightarrow I(y))$, and let $J \equiv I^{\subseteq p}$.

For $QT^+ \vdash J(a)$, note that $QT^+ \vdash I(a)$ since I is a string concept, and $QT \vdash y \subseteq_p a \leftrightarrow y=a$ from (QT2). Hence $QT^+ \vdash \forall y \; (y \subseteq_p a \rightarrow I(y))$. But this suffices for $QT^+ \vdash J(a)$ because $QT \vdash z \leq a \leftrightarrow z=a$.

Likewise $QT^+ \vdash J(b)$, where we need only note that

$$QT \vdash z \leq b \leftrightarrow z=a \lor z=b$$

and appeal to $QT^+ \vdash I(b)$.

Suppose $J(x)$. If $x=a$, we have $Sx=b$, and so $J(Sx)$ by what we just proved.

Otherwise $Sx=x*b$. Suppose $z \leq Sx$, and let $y \subseteq_p z$. It is sufficient to consider the two cases $z \leq x$ and $z = Sx$.

If $z \leq x$, then $I(y)$ follows from the hypothesis $J(x)$.

So let $z = Sx = x*b$. Then, by definition,

$y \subseteq_p x*b \leftrightarrow y = x*b \vee yB(x*b) \vee yE(x*b) \vee \exists x_1, x_2\ x_1 y x_2 = x*b$.

Of the four we only consider the last (see [3], (3.13) for the rest): $x_1 y x_2 = x*b$.

Then $b = x_2 \vee bEx_2$ by (QT4) and (QT5).

Then $x_1 y x_2 = x_1 y b$ or $\exists x_2'\ x_1 y x_2 = x_1 y(x_2'b)$, whence $xb = x_1 y b$ or $xb = x_1 y(x_2'b)$. But then $x = x_1 y$ or $x = x_1 y x_2'$ by (QT3).

In either case $y \subseteq_p x$ and we have $I(y)$ from the hypothesis $J(x)$.

We thus have $M \vDash \forall y\ (y \subseteq_p Sx \rightarrow I(y))$, which is what was needed to show that $\forall z \leq Sx\ \forall y\ (y \subseteq_p z \rightarrow I(y))$.

So we proved that $J(Sx)$ if $J(x)$ as required. That $J(x)$ is closed under $S_a$ is established in a similar fashion. Hence $J(x)$ is indeed a string concept with the required properties.

This means that in establishing that a given string concept I may be strengthened to a string concept J with another property, we need not worry whether the formula J(x) is also closed with respect to * or downward closed with respect to ≤ or ⊆$_p$ . As we just saw, we can always strengthen J(x) to one that is. This is of crucial importance in the formal arguments that we'll be using below, in particular to ensure closure of string concepts under certain kinds of existential claims.

## 6. Adjunctive Set Theory

Let us now consider a very simple set theory, probably the simplest that comes to mind, consisting of the following two principles:

(NULL)  ∃x∀y ¬y∈x,

(ADJ)  ∀x,y∃z∀w (w∈z ↔ (w∈x v w=y))

Formally, we take (NULL) and (ADJ), along with the usual axioms for identity, to determine a first-order theory, Adjunctive Set Theory, AST, formulated in the language $\mathcal{L}$ = {∈}.

By extending AST with

  (EXT) $\qquad\qquad \forall x,y\, (\forall z(z\in x \leftrightarrow z\in y) \rightarrow x=y)$

as an additional axiom, we obtain Adjunctive Set Theory with Extensionality, AST+EXT.

We want to interpret the language $\mathcal{L}$ in a very concrete way. We would like to be able to think of the variables as ranging over non-empty strings of a's and b's, or 0's and 1's. And we want to think of set membership as the relation of one string being a part or a substring of another string: for example, if x=aa and y=baab, then x is a member of y, x∈y, because x is part of y = b |aa |b, or something like that.

## 7. Coding sets by strings

In [13], Quine introduced a method for representing sequences of positive integers by strings of this sort. If we let a tally of n consecutive b's stand for n>0, an ordered pair (i,j) is represented by the string b...bab...b consisting of i many b's followed by a single a and then j many b's. Then a sequence of ordered pairs $(i_1,j_1), (i_2,j_2)...,(i_n,j_n)$ can be represented by the juxtaposition

$$aaw_1aaw_2aa...aaw_naa,$$

of the representations of the pairs separated by aa. This gives us a coding of finite sets of positive integers (or <u>pairs</u> of positive integers) by finite binary strings, made possible by the fact that we can use the a-tallies aa as markers to separate off the members of the set.

But what if $w_1,...,w_n$ were not b-tallies or pairs of b-tallies, but instead arbitrary strings of a's and b's? Here we are facing the problem that no single choice of some string as a marker would suffice because any such string could appear in any of $w_j$'s, or can actually be one of them. So how would we code sets of arbitrary binary strings by single strings?

Suppose t is a b-tally <u>longer than any b-tally occurring</u> in $w_1,...,w_n$.

Interpose copies of t between the members as follows:

$$w = taw_1ataw_2at...aw_nat.$$

Quine made three observations:

(i) t's cannot occur in any $w_j$ because they are longer than any b-tallies contained in them; neither can any occurrences of t straddle any of those shown because they are separated by a's;

(ii) the only segments of the string w which are immediately preceded by ta and immediately followed by at and do not themselves contain any occurrences of t are $w_1,...,w_n$;

(iii) if x is any string not containing any occurrence of t, and if the string taxat occurs as a part of w, then x must be one of $w_1,...,w_n$.

Quine then defines a string x to be a member of the set coded by a string w if the string taxat occurs as a part of w where t is the longest b-tally in w:

$$x \ \varepsilon \ w \equiv \exists t \subseteq_p w \ (MaxT_b(t,w) \ \& \ \neg t \subseteq_p x \ \& \ (taxat) \subseteq_p w).$$

Here we let

$$MaxT_b(t,w) \equiv Tally_b(t) \ \& \ \forall t'(Tally_b(t') \ \& \ t' \subseteq_p w \rightarrow t' \subseteq_p t)$$

where $Tally_b(t) \equiv \forall y \subseteq_p x \ (Digit(y) \rightarrow y=b)$ and $Digit(y) \equiv x=a \lor y=a$.

This means that, in principle, the operation of string concatenation is all we need to be able to express the necessary and sufficient condition for a single binary string to represent a set of strings.

## 8. Step-ladder coding

Note that it is not necessary to use a single b-tally t longer than any b-tally in the member strings to separate them off. It will in fact be more convenient to use different b-tally markers for this purpose: for the string

$$w = t_1 a w_1 a t_2 a w_2 a t_3 \ldots \quad \text{to encode} \quad w_1, w_2, \ldots$$

it will suffice for $t_1$ not to occur in $w_1$, for $t_2$ not to occur in $w_2$, etc., with the additional requirement that the markers $t_1, t_2, t_3, \ldots$ strictly increase in length. The markers $t_1, t_2, t_3, \ldots$ will serve to <u>frame</u> the members $w_1, w_2, \ldots$ .

The coding works like a step-ladder: starting with the b-tally that precedes the first occurrence of the letter a in w, each next longer b-tally is a successive step of the ladder marking off a frame that corresponds to another member of the coded set. (A similar idea was employed by Visser in [19].)

More precisely, let's define when a b-tally t is longer than any b-tally in x:

$$\text{Max}^+\text{T}_b(t,x) \equiv \text{MaxT}_b(t,x) \;\&\; \neg t \subseteq_p x.$$

We then define when a string u is a <u>preframe indexed by</u> t:

$$\text{Pref}(u,t) \equiv \text{Tally}_b(t) \;\&\; \exists y \subseteq_p u \; (aya = u \;\&\; \text{Max}^+\text{T}_b(t,u));$$

when $t_1 u t_2$ is (the) <u>first frame</u> in the string x, $\text{Firstf}(x,t_1,u,t_2)$:

$\text{Pref}(u,t_1) \;\&\; \text{Tally}_b(t_2) \;\&\; ((t_1 = t_2 \;\&\; t_1 u t_2 = x) \lor (t_1 < t_2 \;\&\; (t_1 u t_2 a) B x));$

when $t_1 u t_2$ is (the) <u>last frame</u> in x, $\text{Lastf}(x,t_1,u,t_2)$:

$\text{Pref}(u,t_1) \;\&\; \text{Tally}_b(t_2) \;\&\; t_1 = t_2 \;\&\; (t_1 u t_2 = x \lor \exists w \, (w a t_1 u t_2 = x \;\&\; \text{Max}^+\text{T}_b(t_1,w)));$

and when $t_1 u t_2$ is an <u>intermediate frame</u> in x immediately following an initial segment w of x, $\text{Intf}(x,w,t_1,u,t_2)$:

$\text{Pref}(u,t_1) \;\&\; \text{Tally}_b(t_2) \;\&\; t_1 < t_2 \;\&\; \exists w_1 (w a t_1 u t_2 a w_1 = x) \;\&\; \text{Max}^+\text{T}_b(t_1,w).$

We then define when a string u is $t_1,t_2$-<u>framed</u> in x:

$$\text{Fr}(x,t_1,u,t_2) \equiv \text{Firstf}(x,t_1,u,t_2) \lor \exists w \, \text{Intf}(x,w,t_1,u,t_2) \lor \text{Lastf}(x,t_1,u,t_2),$$

We say that $t_1$ is the <u>initial</u>, and $t_2$ <u>terminal tally marker</u> in the frame.

We then define "t <u>envelops</u> x", Env(t,x), to be the conjunction of the following five conditions:

(a) $\text{MaxT}_b(t,x)$                                "t is a longest b-tally in x",

(b) $\exists u \subseteq_p x \; \exists t_1, t_2 \, \text{Firstf}(x,t_1,u,t_2)$      "x has a first frame",

(c) $\exists u \subseteq_p x$ Lastf(x,t,u,t)          "x has a last frame with t as its initial

                                                                       and terminal marker"

(d) $\forall u \subseteq_p x \; \forall t_1,t_2,t_3,t_4$ (Fr(x,t$_1$,u,t$_2$) & Fr(x,t$_3$,u,t$_4$)) → t$_1$=t$_3$)

      "different initial tally markers frame distinct strings",

(e) $\forall u_1,u_2 \subseteq_p x \; \forall t',t_1,t_2$ (Fr(x,t',u$_1$,t$_1$) & Fr(x,t',u$_2$,t$_2$)) → u$_1$=u$_2$)

      "distinct strings are framed by different initial tally markers"

We then say  x <u>is a set code</u>  if x is aa or else x is enveloped by some b-tally:

   Set(x) ≡ x=aa ∨ ∃t⊆$_p$x Env(t,x).

Finally, we say that a string y <u>is a member of the set coded by</u> string x if x is enveloped by some b-tally t and the juxtaposition of the string y with single tokens of digit a is framed in x:

     y ε x ≡ ∃t⊆$_p$x (Env(t,x) & ∃u⊆$_p$x ∃t$_1$,t$_2$(Fr(x,t$_1$,u,t$_2$) & u=aya)).

Now, suppose a set of strings X is extended by adjoining a string y to obtain another set Y = X ∪ {y}. Then a code for Y can be picked so that a given code of X will be its <u>initial segment</u>.

To be clear, so far we have been talking about codes of sets of strings <u>informally</u>. Various claims were made about properties of codings based on what we took to be obvious properties of concatenated strings. What specific assumptions about strings suffice to formally validate these claims? It will pay off to make these assumptions as weak as possible.

We will show that, modulo our methodology of formula selection, all the necessary reasoning can be carried out in QT$^+$. (See [3], pp.89-263.)

## 9. Interpreting adjunction

We first focus on tallies.

We can show in QT$^+$ that for a suitable string concept J, call it $I_{TOT}$, tallies are totally or completely ordered with respect to $\leq$ ([3], (4.6)):

LEMMA ON COMPLETE ORDERING OF TALLIES. For any string concept $I \subseteq I_0$ there is a string concept $J \subseteq I$ such that

   QT$^+$ ⊢ $\forall z \in J \ \forall x$ (Tally$_b$(x) & Tally$_b$(z) → x$\leq$z v z$\leq$x).

Similarly, we can show that the concatenation operation * is commutative on b-tallies in a suitably defined string concept ([3], (4.10)).

Following this basic methodology, we can obtain progressively more refined string concepts that, provably in our theory, simultaneously have each one from a range of properties needed for our coding of sets to work. A word of caution, however, is in order here. Many seemingly obviously true claims suggested by the formal definitions we gave in the previous section, such as, e.g., that tallies are closed concatenation, the existence and uniqueness of a maximal tally in a given string, the uniqueness of initial and last frames in a given set code, etc., become potentially problematic in the deductively weak setting such as that of QT$^+$. In general, not being provable in unrestricted form in QT$^+$, they must be explicitly proved by selecting appropriate string concepts. (See e.g., [3], (4.5), (4.13) and (4.15).)

But, after numerous auxiliary preparatory such steps it turns out that we can show that any string concept stronger than I$_0$ can be strengthened to one that is closed under set adjunction of strings ([3], (7.1)):

<u>SET ADJUNCTION LEMMA.</u>  For any string concept I ⊆ I$_0$ there is a string concept  J ⊆ I such that
   QT$^+$ ⊢ ∀x,y ∈ J (Set(x) → ∃z ∈ J (Set(z) & ∀w (w ε z ↔ w ε x ∨ w=y))).

On the other hand, from the definitions of the coding scheme we gave earlier,

$$QT^+ \vdash \forall z \,[Set(z) \rightarrow (z=aa \vee \exists y \, y \,\varepsilon\, z) \,\&\, \neg(z=aa \,\&\, \exists y \, y \,\varepsilon\, z)],$$

we can derive that the string aa codes the empty set ([3], (5.18)):

THE NULL SET LEMMA.   $QT^+ \vdash \exists z \,(Set(z) \,\&\, z=aa \,\&\, \forall y \,\neg(y \,\varepsilon\, z)).$

Let the predicate $Set^+(x)$ apply to the set codes among the strings in J,

$$Set^+(x) \equiv J(x) \,\&\, Set(x),$$

where J is obtained from $I_0$ as in the Set Adjunction Lemma.

We define a map $^+$ on atomic formulae of the language of set theory $\mathcal{L} = \{\in\}$ as follows:

$$[x \in y]^+ \equiv Set^+(x) \,\&\, x \,\varepsilon\, y \quad \text{and} \quad [x = y]^+ \equiv x=y.$$

If we let the formula $Set^+(x)$ define the domain of the interpretation, and extend the map $^+$ to non-atomic formulae in the usual way, then the translations of (NULL) and (ADJ)

$$[\exists x \forall y \,\neg y \in x]^+ \quad \text{and} \quad [\forall x,y \exists z \forall w \,(w \in z \leftrightarrow (w \in x \vee w=y))]^+$$

are easily derived ([3], (13.1) and (13.2)):  we have that

$$QT^+ \vdash \exists x \,(Set^+(x) \,\&\, \forall y(Set^+(y) \rightarrow \neg(\, Set^+(y) \,\&\, y \,\varepsilon\, x))).$$

And, from the Set Adjunction Lemma, that

$QT^+ \vdash \forall x,y\ (Set^+(x)\ \&\ Set^+(y) \to \exists z(Set^+(z)\ \&\ \forall w(Set^+(w) \to$

$\to (\ Set^+(w)\ \&\ w\ \varepsilon\ z\ \leftrightarrow\ (Set^+(w)\ \&\ w\ \varepsilon\ x)\ \vee\ w=y))))$.

We thus obtain a formal interpretation of Adjunctive Set Theory (AST) in $QT^+$.

## 10. Canonical set codes

So far we have ignored Extensionality.

It is easy to see that, under the coding scheme we adopted, one and the same set of strings $w_1,...,w_n$ can have different set codes. As we form a string that codes the set, we may take up the members in a different order, or else we may use a different pick of tally markers $t_1,t_2,t_3,...$ to separate the members.

Furthermore, a set code for $w_1,...,w_n$ may contain as substrings material xyz other than the framed strings $w_1,...,w_n$:

$$t_1 a w_1 a t_2 xyz t_3 a w_2 a t_4 ... .$$

This is similar to the so-called "junk DNA" in the human genome.

To validate Extensionality, we'll need to be able to associate a unique string as the canonical set code for strings $w_1,...,w_n$.

Let $Rt_L(z,x,y)$ read "z is the (left) root of x and y":

$(((zaBx \lor za=x) \& (zbBy \lor zb=y)) \lor ((zbBx \lor zb=x) \& (zaBy \lor za=y)))$.

Unless one of the strings x, y, is an initial segment of the other, this says, in effect, z is the longest initial segment common to both x and y. The existence (when it exists) and uniqueness of left root of given strings x and y must be proved by selecting appropriate string concepts (see ([3], (6.2) and (6.3))):

LEFT ROOT LEMMA. For any string concept $I \subseteq I_0$ there is a string concept $I_{RtL} \subseteq I$ such that
$QT^+ \vdash \forall x \in I_{RtL} (\exists z\ zBx \to \forall y\ (x \neq y \to y=a \lor y=b \lor$
$\lor (aBx \& bBy) \lor (bBx \& aBy) \lor xBy \lor yBx \lor \exists z\ Rt_L(z,x,y)))$.

Next, we say that a string u <u>lexically precedes</u> a string v, $u \ll v$, if u is or begins with the letter a and v is or begins with the letter b, or else u is an initial segment of v, or else u and v have a left root z such that u is or begins with za and v is or begins with zb:

$$((u=a \lor aBu) \,\&\, (v=b \lor bBv)) \lor uBv \lor$$
$$\lor \exists z \, (Rt_L(z,u,v) \,\&\, ((za=u \lor zaBu) \,\&\, (zb=v \lor zbBv))).$$

We then have ([3], (6.5)-(6.7)):

LEXICAL PRECEDENCE LEMMA.

(1) For any string concept $I \subseteq I_0$ there is a string concept $J \subseteq I$ such that

$$QT^+ \vdash \forall u,v \in J \, (u \ll v \lor u=v \lor v \ll u).$$

(2) For any string concept $I \subseteq I_0$ there is a string concept $J \subseteq I$ such that

$$QT^+ \vdash \forall v \in J \, \forall u,w \, (u \ll v \,\&\, v \ll w \rightarrow u \ll w).$$

(3) For any string concept $I \subseteq I_0$ there is a string concept $J \subseteq I$ such that

$$QT^+ \vdash \forall u,v \in J \, (u \ll v \rightarrow \neg (v \ll u)).$$

Let's define when a b-tally t is a <u>shortest b-tally not occurring in</u> string u:

$$MinMax^+T_b(t,u) \equiv Max^+T_b(t,u) \,\&\, \forall t'(Max^+T_b(t',u) \rightarrow t \leq t').$$

We can then establish that we can always obtain string concepts in which for every string <u>there does exist</u> a unique shortest b-tally not occurring in that string ([3], (6.8)):

<u>SHORTEST NON-OCCURRENT TALLY LEMMA</u>. For any string concept $I \subseteq I_0$ there is a string concept $J \subseteq I$ such that

$$QT^+ \vdash \forall x \in J \; \exists! t \in J \; MinMax^+T_b(t,x).$$

(We read "$\exists! x \in K \; (\ldots)$" as "$\exists x \; (K(x) \; \& \; (\ldots) \; \& \; \forall y(K(y) \; \& \; (\ldots) \rightarrow y=x))$".)

We now introduce another relation, <u>the shortest b-tally not in u is shorter than the shortest b-tally not in v</u>:

$$u \triangleleft_{Tb} v \equiv \exists t_1, t_2 \; (MinMax^+T_b(t_1,u) \; \& \; MinMax^+T_b(t_2,v) \; \& \; t_1 < t_2).$$

Also, we say when <u>the shortest b-tally not in u is the same as the shortest b-tally not in v</u>:

$$u \approx_{Tb} v \equiv \exists t_1, t_2 \; (MinMax^+T_b(t_1,u) \; \& \; MinMax^+T_b(t_2,v) \; \& \; t_1 = t_2).$$

We can show that any two strings in an appropriately defined string concept are strictly comparable with respect to the shortest b-tallies not occurring in them ([3], (8.1)):

LEMMA ON COMPARABILITY W.R. TO THE SHORTEST NON-OCCURRENT

TALLY. For any string concept $I \subseteq I_0$ there is a string concept $J \subseteq I$ such that

$$QT^+ \vdash \forall u,v \in J((u \triangleleft_{Tb} v \lor u \approx_{Tb} v \lor v \triangleleft_{Tb} u) \& \neg(u \triangleleft_{Tb} v \& v \triangleleft_{Tb} u)).$$

We now order strings accordingly, with the additional stipulation that the strings whose shortest non-occurrent b-tallies are the same are to be ordered according to lexical precedence.

We call this the <u>tally modified lexicographic ordering</u>:

$$u \prec v \equiv (u \triangleleft_{Tb} v \lor (u \approx_{Tb} v \& u \ll v)).$$

We then obtain ([3], (8.2) and (8.3)):

MODIFIED LEXICOGRAPHIC ORDER LEMMA.

(1) For any string concept $I \subseteq I_0$ there is a string concept $J \subseteq I$ such that
$$QT^+ \vdash \forall u,v \in J((u \prec v \lor u = v \lor v \prec u) \& \neg(u \prec v \& v \prec u)),$$

(2) For any string concept $I \subseteq I_0$ there is a string concept $J \subseteq I$ such that
$$QT^+ \vdash \forall v \in J \, \forall u,w \, (u \prec v \& v \prec w \rightarrow u \prec w),$$

Consider now some set code x.

Strings u, v that are members of the set coded by x are embedded in x within frames.

We say that u's frame precedes v's frame in x when either u's frame is the first frame in x, or v's frame is the last frame in x, or else both frames are intermediate and the initial tally marker of v's frame is not shorter than the terminal tally marker of u's frame.

We write  $u<_x v$  for

$\exists t_1,t_2,t_3,t_4[Fr(x,t_1,aua,t_2)\ \&\ Fr(x,t_3,ava,t_4)\ \&$
    $\&\ ((Firstf(x,t_1,aua,t_2)\ \&\ t_1 \neq t_3)\ \vee\ (Lastf(x,t_3,ava,t_4)\ \&\ t_1 \neq t_3)\ \vee$
    $\vee\ (\exists w_1(Intf(x,w_1,t_1,aua,t_2)\ \&\ \exists w_3(Intf(x,w_3,t_3,ava,t_4)\ \&\ t_2 \leq t_3)))].$

We can now state one of the requirements for a set code to count as canonical: the order in which the set members' frames appear in the set's code will have to respect the members' tally modified lexicographic ordering: let

$$Lex^+(z)\ \equiv\ \forall u,v\ (u<_z v \rightarrow u \prec v).$$

Call such set codes lexicographically ordered.

Let's turn now to the analogue of the "no junk DNA" condition, which is a minimality requirement on set codes. We want to make sure that the string coding a given set contains nothing but the framed members of the set. Because tallies of b's serve as markers separating the set's members, the key lies in where we allow the letter a to occur throughout the set code. Given that each frame is of the form

$$t_1 auat_2,$$

we'll let the digit a occur only immediately after an initial tally marker $t_1$, or immediately before a terminal tally marker $t_2$, or else within the framed string u.

We first define when <u>an occurrence of a string z sandwiched between two substrings in x appears within a frame</u>:

x is the result of juxtaposing $w_1$ to the left and $w_2$ to the right of z, and either $t_1 v t_2$ is the first frame in x, and the string $w_1$ is the b-tally $t_1$ or $t_1 v$ begins with or is the string $w_1 z$; or else $t_1 v t_2$ is an intermediate or the last frame in x having some string w'a as the initial segment of x immediately preceding it, and the string $w_1$ is w'a$t_1$, or $w_1 z$ is w'a$t_1 v$ or results from juxtaposing some initial segment $v_1$ of v next to w'a $t_1$.

That is, we write  $Occ(w_1,z,w_2,x,t_1,v,t_2)$  for

$\quad w_1zw_2=x$ & $Fr(x, t_1,v,t_2)$ &

& $[(Firstf(x, t_1,v,t_2)$ & $(t_1=w_1 \lor (w_1z)B(t_1v) \lor w_1z= t_1v)) \lor$

$\lor \exists w'((Intf(x,w', t_1,v,t_2) \lor (Lastf(x, t_1,v,t_2)$ & $w'at_1vt_2=x))$ &

& $(w'at_1=w_1 \lor \exists v_1(v_1Bv$ & $w'a\, t_1v_1=w_1z) \lor w'at_1v=w_1z))]$.

We call a set code <u>minimal</u> if every occurrence of the digit a appears within some frame: write  $MinSet(x)$  for

$Set(x)$ & $\forall w_1,w_2 \subseteq_p x\ (w_1aw_2=x \rightarrow \exists v \subseteq_p x\ \exists t',t'' \subseteq_p x\ Occ(w_1,a,w_2,x,t',v,t''))$.

Let's recall how our step-ladder coding works.

Given strings $w_1,w_2,\ldots$ we select tally markers $t_1,t_2,t_3,\ldots$ for the corresponding frames to obtain

$$w = t_1aw_1at_2aw_2at_3\ldots.$$

The tally markers are strictly increasing in length: we have to make sure that the initial tally marker for a string's frame is longer than any b-tally in that frame, and also longer than any initial tally marker corresponding to frames that precede that string's frame in the code.

But this leaves us free to make the initial tally markers as large as we want.

For canonical set codes we require that the b-tallies used as initial markers be shortest possible.

For a string v in a given set of strings to be coded, first we state the condition a b-tally t must meet to serve as a possible initial tally marker for v's frame in some set code x: the b-tally t should be longer than the initial tally marker of any frame that precedes v's frame in x.

Write $Max^+(t,v,x)$ for

$Set(x)$ & $v \; \varepsilon \; x$ & $Tally_b(t)$ &
    & $\forall u, t_1, t_2 \subseteq_p x \; (Fr(x, t_1, aua, t_2) \; \& \; u <_x v \rightarrow t_1 < t)$.

We then require that t be the shortest such tally:

Let $MMax^+T_b(t,v,x) \equiv Max^+(t,v,x) \; \& \; \forall t'(Max^+(t',v,x) \; \& \; Max^+T_b(t',v) \rightarrow t \leq t')$.

We call set codes in which each frame has as its initial tally marker a b-tally that (uniquely) satisfies this condition special.

Let $Special(x) \equiv Set(x) \; \& \; \forall v, t_1, t_2 \; (Fr(x, t_1, ava, t_2) \rightarrow MMax^+T_b(t_1,v,x))$.

If we let ~ mean that the sets coded by x and y have the same strings as members,

$$x \sim y \equiv Set(x) \ \& \ Set(y) \ \& \ \forall w(w \ \varepsilon \ x \leftrightarrow w \ \varepsilon \ y),$$

we then have the Special Set Codes Lemma, which says that we can choose a string concept in which, for set codes that are both lexicographically ordered and special, the members of the coded sets uniquely determine the initial tally markers of their frames ([3], (10.2)):

SPECIAL SET CODES LEMMA.  For any string concept $I \subseteq I_0$ there is a string concept $J \subseteq I$ such that
$$QT^+ \vdash \forall x,y \in J \ (Lex^+(x) \ \& \ Lex^+(y) \ \& \ Special(x) \ \& \ Special(y) \ \& \ x \sim y \ \&$$
$$\& \ Fr(x,t_1,aua,t_2) \ \& \ Fr(y,t_3,aua,t_4) \rightarrow t_1 = t_3).$$

We call a set code <u>canonical</u> if its lexicographically ordered, minimal and special:

$$Set^*(x) \equiv MinSet(x) \ \& \ Lex^+(x) \ \& \ Special(x).$$

We then have ([3], (11.4)):

THE UNIQUENESS LEMMA. For any string concept $I \subseteq I_0$ there is a string concept $J \subseteq I$ such that

$$QT^+ \vdash \forall x,y \in J \; (Set^*(x) \; \& \; Set^*(y) \; \& \; x \sim y \; \rightarrow \; x=y).$$

Sets that have the same strings as members have the same canonical set code.

## 11. Interpreting Extensionality

The Uniqueness Lemma will be an essential element in our interpretation of (EXT). But we simultaneously have to make sure that (ADJ) also holds:

STRONG SET ADJUNCTION LEMMA. There is a formula $V^{**} \subseteq Set^*$ such that
$$QT^+ \vdash \forall x,y \; (V^{**}(x) \; \& \; V^{**}(y) \; \rightarrow \; \exists!z \; (V^{**}(z) \; \& \; \forall w \; (w \; \varepsilon \; z \leftrightarrow (w \; \varepsilon \; x \; \vee \; w=y)))).$$

In contrast to the version of the Set Adjunction Lemma we used earlier to interpret AST without EXT, here the canonical code $z$ produced by adjunction cannot be obtained simply by tacking on a frame for the adjoined string at the tail end of the set code x for the original set. The canonical code for the expanded set has to be reconfigured using a suitable selection of initial tally markers, depending on how the new member y lexicographically relates to the members of the original set. The proof of the Strong Set Adjunction Lemma requires that we consider the whole variety of cases that arise in this connection (see [3], (12.2)). From the proof we can extract a rather lengthy formula $\sigma^*(x,y,z)$ in the language of $QT^+$ for which we obtain (see [3], pp.646-650):

STRONG SET ADJUNCTION LEMMA (EXPLICIT FORM).

$QT^+ \vdash \forall x,y\ (V^{**}(x)\ \&\ V^{**}(y) \rightarrow \exists!z\ (V^{**}(z)\ \&\ \sigma^*(x,y,z))\ \&$

$\&\ \forall z(\sigma^*(x,y,z) \rightarrow \forall w\ (w\ \varepsilon\ z \leftrightarrow (w\ \varepsilon\ x \lor w=y))\ \&\ (\sigma^*(x,y,x) \leftrightarrow y\ \varepsilon\ x)).$

Now we are ready to set up our formal interpretation.

We'll let the formula $V^{**}(x)$ define the domain.

We write $x \stackrel{*}{=} y$, <u>x and y code the same set modulo $V^{**}$</u>, for

$$\forall z(V^{**}(z) \rightarrow (z\ \varepsilon\ x \leftrightarrow z\ \varepsilon\ y)).$$

We interpret atomic formulae of the language of set theory $\mathcal{L} = \{\in\}$ as follows:

$[x=y]^* \equiv x \stackrel{*}{=} y$ and $[x \in y] \equiv x\ \varepsilon\ y.$

We then extend the map * to non-atomic formulae in the usual way,

relativizing the quantifiers to V**.

Then the *-translation of (NULL), [∃x∀y ¬y∈x]*, is proved in the same way as the +-translation [∃x∀y ¬y∈x]+ earlier.

On the other hand, from the Strong Set Adjunction Lemma we show that

QT+ ⊢ ∀x (V**(x) → ∀y(V**(y) → ∃z(V**(z) & y ε z &
           & ∀w(V**(w) → (w ε x → w ε z)) &
    & ∀w(V**(w) → (w ε z → w ε x ∨ ∀v(V**(v) → (v ε w ↔ v ε y)))))))).

But this is the *-translation of a formula equivalent to the Adjunction axiom:

QT+ ⊢ [∀x,y∃z (y ε z & ∀w (w∈x → w∈z) & ∀w (w∈z → w∈x ∨ w=y))]*.

Finally, note that

∀x (V**(x) → ∀y(V**(y) → (∀z(V**(z) → (z ε x ↔ z ε y)) →
                     → ∀z(V**(z) → (z ε x ↔ z ε y)))))

is in fact the *-translation [∀x,y (∀z(z∈x ↔ z∈y) → x=y)]* of the Extensionality axiom.

Hence QT+ ⊢ [EXT]* holds trivially.

So our concatenation theory QT+ also interprets Adjunctive Set Theory with Extensionality, AST+EXT (see [3], pp.662-663):

Theorem 1.    AST+EXT $\leq_I$ QT$^+$.

## 12. Concatenation arithmetically represented

The arithmetical theory Q, known as Robinson Arithmetic, is formulated in the first-order language $\{0, ', +, \cdot\}$, with (the universal closures of) the following (non-logical) axioms:

(Q1)   $\neg x'=0$

(Q2)   $x'=y' \rightarrow x=y$

(Q3)   $x=0 \lor \exists y\, y'=x$

(Q4)   $x+0=x$

(Q5)   $x+y'=(x+y)'$

(Q6)   $x \cdot 0=0$

(Q7)   $x \cdot y'=x \cdot y+x$

Now, QT$^+$ can be interpreted in the arithmetical theory

$$I\Sigma_0 = Q + \{\, `\varphi(0) \,\&\, \forall x\,(\varphi(x) \to \varphi(Sx)) \to \forall x\,\varphi(x)\text{'} \mid \varphi(x) \in \Sigma_0 \,\}$$

where $\varphi(x)$ is any bounded formula in the same language $\{0,',+,\cdot\}$, i.e. formula with no unbounded quantifiers, but possibly with parameters other than x.

In $I\Sigma_0$ it is possible to define a coding of sequences of numbers and a corresponding concatenation operation: specifically, there is

(i) a bounded formula, *Seq*(x), that defines the set of numbers that serve as the codes of sequences of numbers, including a code for the empty sequence,

(ii) a bounded predicate, $x \in y$, expressing the relation "x is a term of the sequence coded by y",

(iii) a polynomially bounded function $x \frown y$ that yields the code of the concatenation of two sequences, the sequence whose terms are the terms of the sequence coded by x followed by the terms of the sequence coded by y, provided *Seq*(x) and *Seq*(y) (otherwise, $x \frown y = 0$).

We let $Seq^*(x) \equiv Seq(x) \,\&\, \exists y\, y \in x \,\&\, \forall y\,(y \in x \to y = S0 \,\lor\, y = SS0)$,

Then the predicate *Seq*\*(x) defines the set of codes of nonempty dyadic sequences, i.e. sequences of 1's and/or 2's. That is, we have:

(t0)  $I\Sigma_0 \vdash \exists s\, Seq^*(s)$

(t1)  $I\Sigma_0 \vdash \exists!s\ (Seq^*(s)\ \&\ \forall x\ (x \in s \leftrightarrow x = \underline{c_1}))$

(t2)  $I\Sigma_0 \vdash \exists!s\ (Seq^*(s)\ \&\ \forall x\ (x \in s \leftrightarrow x = \underline{c_2}))$

(t3)  $I\Sigma_0 \vdash Seq^*(s)\ \&\ Seq^*(t) \rightarrow \exists!u\ (Seq^*(u)\ \&\ s\frown t = u)$

(t4)  $I\Sigma_0 \vdash Seq^*(s)\ \&\ Seq^*(t)\ \&\ Seq^*(u) \rightarrow (s\frown t)\frown u = s\frown(t\frown u)$

(t5)  $I\Sigma_0 \vdash Seq^*(s)\ \&\ Seq^*(t) \rightarrow \neg(s\frown t = \underline{c_1})\ \&\ \neg(s\frown t = \underline{c_2})$

(t6)  $I\Sigma_0 \vdash Seq^*(s)\ \&\ Seq^*(t) \rightarrow (\underline{c_1}\frown s = \underline{c_1}\frown t \rightarrow s = t)\ \&\ (\underline{c_2}\frown s = \underline{c_2}\frown t \rightarrow s = t)$

(t7)  $I\Sigma_0 \vdash Seq^*(s)\ \&\ Seq^*(t) \rightarrow (s\frown \underline{c_1} = t\frown \underline{c_1} \rightarrow s = t)\ \&\ (s\frown \underline{c_2} = t\frown \underline{c_2} \rightarrow s = t)$

(t8)  $I\Sigma_0 \vdash Seq^*(s)\ \&\ Seq^*(t) \rightarrow \neg(\underline{c_1}\frown s = \underline{c_2}\frown t)\ \&\ \neg(s\frown \underline{c_1} = t\frown \underline{c_2})$

(t9)  $I\Sigma_0 \vdash Seq^*(s) \rightarrow s = \underline{c_1}\ \vee\ s = \underline{c_2}\ \vee\ (\exists t\ (Seq^*(t)\ \&\ (\underline{c_1}\frown t = s\ \vee\ \underline{c_2}\frown t = s))\ \&$
$\&\ \exists t\ (Seq^*(t)\ \&\ (t\frown \underline{c_1} = s\ \vee\ t\frown \underline{c_2} = s)))$

Here $\underline{c_1}$ and $\underline{c_2}$ are variable free terms that are the codes of the singleton sequences that consist of 1 and 2, respectively. In particular, (t4)-(t9) verify the axioms of QT in $I\Sigma_0$. This determines an interpretation of concatenation theory QT in $I\Sigma_0$.

Since QT$^+$ is interpretable in QT, this also establishes

Theorem 2.    $QT^+ \leq_I I\Sigma_0$.

But $I\Sigma_0$ is known to be interpretable in Q, by a theorem of Wilkie (see [8]). Hence we also have

Theorem 3.    $QT^+ \leq_I Q$.

## 13. All the pieces fall into place

By Tarski and Szmielew, Collins and Halpern, and Mycielski, Pudlák and Stern (who dispensed with extensionality), Q is relatively interpretable in Adjunctive Set Theory AST. So $Q \leq_I AST$, whereas in Theorem 1 we have shown that $AST+EXT \leq_I QT^+$.

This closes the circle: it follows that Q, $QT^+$, AST and AST+EXT are all mutually interpretable:

FIRST MUTUAL INTERPRETABILITY THEOREM.

$$Q \equiv_I QT^+ \equiv_I AST \equiv_I AST+EXT.$$

## 14. Adjunction in functional form: quantifier-free finitary set theory

The theory $PS_0$, introduced by Kirby in [9] as the quantifier-free theory of finite sets, is formulated in the first-order language {0, ;}, where "0" is an individual constant and ";" a binary function symbol.

Its axioms are

(PS1)  $0;x \neq 0$

(PS2)  $x;y,y = x;y$

(PS3)  $x;y,z = x;z,y$

(PS4)  $x;y,z = x;y \leftrightarrow x;z = x \lor z = y$

where we write "x;y,z" for "(x;y);z".

If we let $\quad\quad\quad x \in y \equiv x;y=x,$

we may informally express the meaning of (PS1)-(PS4) in more familiar notation as:

$\quad \neg(x \in \emptyset)$

$\quad y \in x \cup \{y\}$

$\quad (x \cup \{y\}) \cup \{z\} = (x \cup \{z\}) \cup \{y\}$

$\quad z \in x \cup \{y\} \leftrightarrow z \in x \lor z = y.$

$PS_0$ may be thought of as a minimal theory of the <u>adjunction operation</u> $x \cup \{y\}$.

Note that

$\quad$ (PS2) & (PS3) $\vdash x;z = x \lor z = y \rightarrow x;y,z = x;y.$

For, assume $z = y$. Then $x;y,z = x;y,y = x;y$ by (PS2).

Hence (PS2) $\vdash z = y \rightarrow x;y,z = x;y.$

Assume $x;z = x$. Then $x;z,y = (x;z);y = x;y$. But $x;z,y = x;y,z$ by (PS3).

So $x;y,z = x;y.$

Hence (PS3) ⊢ x;z = x → x;y,z = x;y.

Therefore

   (PS2) & (PS3) ⊢ x;z = x ∨ z = y → x;y,z = x;y, as claimed.

In more familiar notation, (PS2) & (PS3) ⊢ z ∈ x ∨ z = y → z ∈ x ∪ {y}.

Minimal Predicative Set Theory N studied by Montagna and Mancini in [11] amounts to taking (PS1) and (PS4) as sole non-logical axioms along with the usual axioms for identity.

For our purpose it will be convenient to reformulate $PS_0$ as a first-order theory $PS_0'$ in the language {0, S} where "0" is an individual constant and "S" a ternary relation symbol satisfying the following conditions:

(PS1′) ¬S(0,x,0)

(PS2′) $S(x,y,z_1)$ & $S(z_1,y,z_2)$ → $z_1 = z_2$

(PS3′) $S(x,y,z_1)$ & $S(z_1,z,z_2)$ & $S(x,z,z_3)$ & $S(z_3,y,z_4)$ → $z_2 = z_4$

(PS4′) S(x,y,z) & S(z,w,z) → S(x,w,x) ∨ w=y

(PS5′) ∃z S(x,y,z)

(PS6′) $S(x,y,z_1)$ & $S(x,y,z_2)$ → $z_1 = z_2$.

Again, we take the formula V**(x) from the STRONG SET ADJUNCTION LEMMA to define the domain of the interpretation.

We let the string aa interpret the constant 0, and let

$$[S(x,y,z)]^* \equiv \sigma^*(x,y,z),$$

where $\sigma^*(x,y,z)$ is as in the EXPLICIT FORM of the STRONG SET ADJUNCTION LEMMA.

Let $[x = y]^* \equiv x=y$.

Then we can verify:

(i) $QT^+ \vdash V^{**}(x) \to \neg[S(0,x,0)]^*$.

Assume $V^{**}(x)$. Suppose, for a reductio, that $\sigma^*(aa,x,aa)$. By the STRONG SET ADJUNCTION LEMMA,

$$\forall w(w \; \varepsilon \; aa \leftrightarrow w \; \varepsilon \; aa \lor w=x),$$

whence $x \; \varepsilon \; aa \leftrightarrow x \; \varepsilon \; aa \lor x=x$.

But, $QT^+ \vdash \neg(x \; \varepsilon \; aa)$, as noted in the NULL SET LEMMA.

Since $x=x \lor x \; \varepsilon \; aa$, this is a contradiction.

Therefore, $\neg\sigma^*(aa,x,aa)$, that is, $M \vDash \neg[S(0,x,0)]^*$.

Hence $V^{**}(x) \to \neg\sigma^*(aa,x,aa)$.

(ii) $QT^+ \vdash V^{**}(x) \; \& \; V^{**}(y) \; \& \; V^{**}(z_1) \; \& \; V^{**}(z_2) \to$

$$\to [S(x,y,z_1) \; \& \; S(z_1,y,z_2) \to z_1=z_2]^*.$$

Assume $V^{**}(x) \; \& \; V^{**}(y) \; \& \; V^{**}(z_1) \; \& \; V^{**}(z_2)$ along with

$$\sigma^*(x,y,z_1) \; \& \; \sigma^*(z_1,y,z_2).$$

(iia) $M \vDash y \; \varepsilon \; x$.

Then from $\sigma^*(x,y,z_1)$, $x=z_1$. Then $y \; \varepsilon \; z_1$, whence $z_1=z_2$ from $\sigma^*(z_1,y,z_2)$.

(iib) $M \vDash \neg(y \; \varepsilon \; x)$.

From the hypothesis, according to the EXPLICIT FORM of the STRONG SET ADJUNCTION LEMMA we have that $y \varepsilon z_1$. But then $z_1=z_2$ from $\sigma^*(z_1,y,z_2)$. Hence we have

$\quad V^{**}(x) \& V^{**}(y) \& V^{**}(z_1) \& V^{**}(z_2) \rightarrow (\sigma^*(x,y,z_1) \& \sigma^*(z_1,y,z_2) \rightarrow z_1=z_2)$,

that is,

$QT^+ \vdash V^{**}(x) \& V^{**}(y) \& V^{**}(z_1) \& V^{**}(z_2) \rightarrow [S(x,y,z_1) \& S(z_1,y,z_2) \rightarrow z_1=z_2]^*$.

(iii) $QT^+ \vdash V^{**}(x) \& V^{**}(y) \& V^{**}(z) \& V^{**}(z_1) \& V^{**}(z_2) \& V^{**}(z_3) \& V^{**}(z_4) \rightarrow$

$\quad\quad\quad\quad\quad \rightarrow [S(x,y,z_1) \& S(z_1,z,z_2) \& S(x,z,z_3) \& S(z_3,y,z_4) \rightarrow z_2=z_4]^*$.

Assume that $V^{**}(x) \& V^{**}(y) \& V^{**}(z_1) \& V^{**}(z_2) \& V^{**}(z_3) \& V^{**}(z_4)$

and $\quad\quad\quad \sigma^*(x,y,z_1) \& \sigma^*(z_1,z,z_2) \& \sigma^*(x,z,z_3) \& \sigma^*(z_3,y,z_4)$.

We then have that

$\quad \forall w(w \varepsilon z_2 \leftrightarrow w \varepsilon z_1 \vee w=z \leftrightarrow (w \varepsilon x \vee w=y) \vee w=z \leftrightarrow$

$\quad\quad\quad\quad \leftrightarrow (w \varepsilon x \vee w=z) \vee w=y \leftrightarrow w \varepsilon z_3 \vee w=y \leftrightarrow w \varepsilon z_4)$.

That is, $z_2 \sim z_4$.

But from $V^{**}(z_2) \& V^{**}(z_4)$ we have $Set^*(z_2) \& Set^*(z_4)$.

Hence, by the UNIQUENESS LEMMA, $z_2=z_4$,

which suffices to prove the claim.

(iv) $QT^+ \vdash V^{**}(x) \& V^{**}(y) \& V^{**}(z) \& V^{**}(w) \rightarrow$

$\quad\quad\quad\quad\quad \rightarrow [S(x,y,z) \& S(z,w,z) \rightarrow S(x,w,x) \vee w=y]^*$.

Assume  $V^{**}(x)$ & $V^{**}(y)$ & $V^{**}(z)$ & $V^{**}(w)$

along with  $\sigma^*(x,y,z)$ & $\sigma^*(z,w,z)$.

From the STRONG SET ADJUNCTION LEMMA we have

$\forall v(v \, \varepsilon \, z \leftrightarrow w \, \varepsilon \, x \lor v=y)$  and  $\forall v(v \, \varepsilon \, z \leftrightarrow v \, \varepsilon \, z \lor v=w)$.

Also,  $w \, \varepsilon \, z$.

Assume that  $M \vDash w \neq y$.

Then  $w \, \varepsilon \, x$. But then  $\sigma^*(x,w,x)$, as required.

(v) $QT^+ \vdash V^{**}(x)$ & $V^{**}(y) \rightarrow \exists z(V^{**}(z)$ & $[S(x,y,z)]^*)$.

(vi) $QT^+ \vdash V^{**}(x)$ & $V^{**}(y)$ & $V^{**}(z_1)$ & $V^{**}(z_2) \rightarrow$

$$\rightarrow [S(x,y,z_1) \, \& \, S(x,y,z_2) \rightarrow z_1=z_2]^*.$$

This follows immediately from the STRONG SET ADJUNCTION LEMMA.

With (i)-(vi) we have derived:    $PS_0 \leq_I QT^+$.

Let us now consider extensionality in this setting.

The axioms of $PS_0 + EXT$ are those of $PS_0$ together with

  EXT.    $\forall x,y(\forall z(S(x,z,x) \leftrightarrow S(y,z,y)) \rightarrow x=y)$.

Again, let $V^{**}(x)$ define the domain of the interpretation **, where

  $[0]^{**} \equiv aa$  and  $[S(x,y,z)]^{**} \equiv \sigma^*(x,y,z)$.

Let $[x=y]^{**} \equiv x \stackrel{*}{=} y$

where $\quad x \stackrel{*}{=} y \equiv \forall z(V^{**}(z) \to (z \,\varepsilon\, x \leftrightarrow z \,\varepsilon\, y))$.

We then argue just as in (i)-(ii) that

(i**) $QT^+ \vdash V^{**}(x) \to \neg[S(0,x,0)]^{**}$,

(ii**) $QT^+ \vdash V^{**}(x) \& V^{**}(y) \& V^{**}(z_1) \& V^{**}(z_2) \to$

$$\to [S(x,y,z_1) \& S(z_1,y,z_2) \to z_1 = z_2]^{**},$$

(iii***)

$QT^+ \vdash V^{**}(x) \& V^{**}(y) \& V^{**}(z) \& V^{**}(z_1) \& V^{**}(z_2) \& V^{**}(z_3) \& V^{**}(z_4) \to$

$$\to [S(x,y,z_1) \& S(z_1,z,z_2) \& S(x,z,z_3) \& S(z_3,y,z_4) \to z_2 = z_4]^{**}.$$

We argue as in (iii) and obtain, under the hypothesis, that

$$\forall w(w \,\varepsilon\, z_2 \leftrightarrow w \,\varepsilon\, z_4).$$

But then $z_2 \stackrel{*}{=} z_4$.

(iv**) $QT^+ \vdash V^{**}(x) \& V^{**}(y) \& V^{**}(z) \& V^{**}(w) \to$

$$\to [S(x,y,z) \& S(z,w,z) \to S(x,w,x) \lor w = y]^{**}.$$

Assume $V^{**}(x) \& V^{**}(y) \& V^{**}(z) \& V^{**}(w)$ along with $\sigma^*(x,y,z) \& \sigma^*(z,w,z)$.

From the STRONG SET ADJUNCTION LEMMA we have

$$\forall v(v \,\varepsilon\, z \leftrightarrow w \,\varepsilon\, x \lor v = y) \quad \text{and} \quad \forall v(v \,\varepsilon\, z \leftrightarrow v \,\varepsilon\, z \lor v = w).$$

But $w \,\varepsilon\, z$ from the hypothesis. Hence $w \,\varepsilon\, x \lor w = y$.

It follows that $\sigma^*(x,w,x) \lor \forall v(V^{**}(v) \to (v \,\varepsilon\, w \leftrightarrow v \,\varepsilon\, y))$, as required.

(v**) $QT^+ \vdash V^{**}(x) \& V^{**}(y) \to \exists z(V^{**}(z) \& [S(x,y,z)]^{**})$ is

immediate from the STRONG SET ADJUNCTION LEMMA.

(vi**) $QT^+ \vdash V^{**}(x) \& V^{**}(y) \& V^{**}(z_1) \& V^{**}(z_2) \to$

$$\to [S(x,y,z_1) \& S(x,y,z_2) \to z_1=z_2]^{**}.$$

Assume $V^{**}(x) \& V^{**}(y) \& V^{**}(z_1) \& V^{**}(z_2)$ along with $\sigma^*(x,y,z_1) \& \sigma^*(x,y,z_2)$.
From the STRONG SET ADJUNCTION LEMMA we have that

$\forall w(w \; \varepsilon \; z_1 \leftrightarrow w \; \varepsilon \; x \lor w=y)$ and $\forall w(w \; \varepsilon \; z_2 \leftrightarrow w \; \varepsilon \; x \lor w=y)$,

So $\forall w(w \; \varepsilon \; z_1 \leftrightarrow w \; \varepsilon \; z_2)$, whence $\forall w(V^{**}(w) \to (w \; \varepsilon \; z_1 \leftrightarrow w \; \varepsilon \; z_2))$,

as required.

Finally, we have that

(EXT**) $QT^+ \vdash V^{**}(x) \& V^{**}(y) \to [\forall z(S(x,z,x) \leftrightarrow S(y,z,y)) \to x=y]^{**}$.

Assume $V^{**}(x) \& V^{**}(y)$ along with $\forall z(V^{**}(z) \to (\sigma^*(x,z,x) \leftrightarrow \sigma^*(y,z,y)))$.

From the explicit form of the STRONG SET ADJUNCTION LEMMA we have that

$\forall w(V^{**}(w) \to (\sigma^*(x,w,x) \leftrightarrow w \; \varepsilon \; x))$

and also $\forall w(V^{**}(w) \to (\sigma^*(y,w,y) \leftrightarrow w \; \varepsilon \; y))$.

Hence $M \vDash \forall z(V^{**}(z) \to (z \; \varepsilon \; x \leftrightarrow z \; \varepsilon \; y))$, as required.

In deriving (i**)-(vii*) and (EXT**) we have also established the interpretability of $PS_0$ + EXT in $QT^+$.

So we have that $\quad N \leq_I PS_0 \leq_I PS_0 + EXT \leq_I QT^+$.

From Montagna and Mancini [11] we have $Q \leq_I N$.

Since, by Theorem 3 above, $QT^+ \leq_I Q$, the circle closes again and we have

SECOND MUTUAL INTERPRETABILITY THEOREM.

$$Q \equiv_I QT^+ \equiv_I N \equiv_I PS_0 \equiv_I PS_0 + EXT.$$